\documentclass{amsart}

\usepackage[english]{babel}
\usepackage{amsmath,amssymb,amsthm}

\usepackage[T1]{fontenc}
\usepackage[utf8x]{inputenc}
\usepackage{enumitem}
\usepackage{verbatim}
\usepackage{graphicx}
\usepackage{faktor}
\usepackage{xcolor}
\usepackage{xfrac}
\usepackage{hyperref}
\usepackage[normalem]{ulem}
\usepackage{setspace}
\usepackage{bbm}
\usepackage{calrsfs}
\DeclareMathAlphabet{\pazocal}{OMS}{zplm}{m}{n}

\newcommand{\etalchar}[1]{$^{#1}$}

\newcommand{\PP}{\mathbb P}
\newcommand{\CC}{\mathbb C}
\newcommand{\Aaff}{\mathbb A}

\newcommand{\ZZ}{\mathbb Z}

\newcommand{\frakc}{\mathfrak c}
\newcommand{\m}{\mathfrak m}
\newcommand{\tR}{\widetilde{R}}
\newcommand{\tm}{\widetilde{\mathfrak m}}

\newcommand{\OO}{\mathcal O}

\newcommand{\Ccal}{\mathcal C}
\newcommand{\Hcal}{\pazocal H}
\newcommand{\Mcal}{\pazocal M}

\newcommand{\oM}{\overline{\pazocal M}}

\newcommand{\oC}{\overline{C}}
\newcommand{\de}{\operatorname{d}}

\renewcommand{\to}{\rightarrow}

\newcommand{\Gm}{\mathbb{G}_{\rm{m}}}
\newcommand{\Ga}{\mathbb{G}_{\rm{a}}}

\newcommand{\dvr}{\Delta}

\newcommand{\bq}{\begin{equation}}
\newcommand{\eq}{\end{equation}}
\newcommand{\ba}{\begin{aligned}}
\newcommand{\ea}{\end{aligned}}
\newcommand{\be}{\begin{enumerate}}
\newcommand{\ee}{\end{enumerate}}
\newcommand{\bsm}{\left(\begin{smallmatrix}}
\newcommand{\esm}{\end{smallmatrix}\right)}                   
\newcommand{\bpm}{\begin{pmatrix}}
\newcommand{\epm}{\end{pmatrix}}
\newcommand{\barr}{\begin{displaymath}\begin{array}{cccc}}
\newcommand{\earr}{\end{array}\end{displaymath}}
\newcommand{\barrl}{\begin{displaymath}\begin{array}{lcl}}
\newcommand{\earrl}{\end{array}\end{displaymath}}
\newcommand{\barl}{\begin{displaymath}\begin{array}{l}}
\newcommand{\earl}{\end{array}\end{displaymath}}
\newcommand{\bxym}{ \begin{displaymath}\xymatrix }
\newcommand{\exym}{\end{displaymath}}
\newcommand{\bcd}{\begin{center}\begin{tikzcd}}
\newcommand{\ecd}{\end{tikzcd}\end{center}}
\newcommand{\on}{\operatorname}
\newcommand{\cond}{\mathfrak{c}}

\newcommand{\Pic}{\operatorname{Pic}}

\theoremstyle{plain}
\newtheorem{thm}{Theorem}[section]
\newtheorem{lem}[thm]{Lemma}

\newtheorem{conj}[thm]{Conjecture}
\newtheorem*{theorem*}{Theorem}

\theoremstyle{definition}

\theoremstyle{remark}
\newtheorem{remark}[thm]{Remark}
\newtheorem{rem}[thm]{Remark}
\newtheorem{rmk}[thm]{Remark}

\usepackage{tikz,tikz-cd,tikz-3dplot}
\usepackage{pgfplots}
\usetikzlibrary{arrows,shadows,positioning, calc, decorations.markings, 
hobby,quotes,angles,decorations.pathreplacing,intersections,shapes}
\usepgflibrary{shapes.geometric}
\usetikzlibrary{fillbetween,backgrounds}

\usepackage{orcidlink}

\title[Singularities and differentials of genus three]{Gorenstein curve singularities of genus three}
\author[L. Battistella]{Luca Battistella\,\orcidlink{0000-0002-0042-970X}}
\email{luca.battistella2@unibo.it}
\address{Dipartimento di Matematica, Universit\`a di Bologna, Italy}

\thanks{\emph{Acknowledgements.} I would like to thank D. Agostini, S. Bozlee, D. Chen, R. Krishnamoorthy, D. Ranganathan, and J. Wise for helpful conversations, and especially F. Bernasconi, F. Carocci, F. Viviani, and the anonymous referee for precise comments on the manuscript.}


\begin{document}

\begin{abstract}
We classify the analytic germs of isolated Gorenstein curve singularities of genus three, and relate them to the connected components of strata of abelian differentials. %The motivation comes from moduli theory.
\end{abstract}

\maketitle 

\tableofcontents

\section*{Introduction}
A Cohen-Macaulay variety is Gorenstein if its dualising sheaf is a line bundle. This is a local condition. For instance, ordinary double points are Gorenstein, whereas the union of the coordinate axes in $\Aaff^3$ is not. The class of Gorenstein curves provides a very mild and useful generalisation of nodal curves from the moduli theoretic perspective.
\begin{itemize}[leftmargin=0.5cm]
 \item There are several compactifications of $\pazocal M_{g,n}$ parametrising only Gorenstein curves: the most notable one is the Deligne--Mumford compactification $\oM_{g,n}$, and the alternatives serve to illustrate its birational geometry, see for instance \cite{HassettHyeonContr,SMY1,AFS1,Bat19,BKN}. Indeed, all singularities appearing in the Hassett--Keel program seem to be Gorenstein \cite{AFSGm}. Typically, replacing nodal subcurves with few special points with a Gorenstein singularity produces a contraction or flip of $\oM_{g,n}$, thus worsening the singularities of the moduli space, but otherwise reducing its complexity (the number of parameters).
 \item Gorenstein curves can be used to test the geometry of a target variety \cite{Viscardi,BCM}. %The deformation theory of maps is identical to the nodal case \cite[Proposition 6.2]{BF}, resulting into a virtual fundamental class.
 They have been successfully employed to construct alternative compactifications and modular desingularisations of Kontsevich's moduli space of stable maps \cite{RSPW1,RSPW2,BNR1,BCquartics,BatCar}, thus defining simpler enumerative invariants.
\end{itemize}

Let us now consider the analytic germ of an isolated curve singularity $(C,p)$. The simplest of its invariants are the following:
\begin{itemize}[leftmargin=0.5cm]
 \item $b$, the number of branches (or preimages of $p$ in the normalisation);
 \item $\delta$, the number of independent conditions that a function on the normalisation must satisfy in order to descend to the singular curve. Out of these, $b-1$ impose that the function attains the same value at the $b$ preimages of $p$;
 \item the remaining $g:=\delta-b+1$ form the local contribution of the singularity to the arithmetic genus of any projective curve containing it; $g$ is called the \emph{genus} of the singularity.
\end{itemize}

Several more sophisticated invariants of curve singularities have been introduced, see for instance the book \cite{Greuel}, although most classical sources are only concerned with the singularities of curves on surfaces, see e.g. the book \cite{Brieskorn}. From the moduli-theoretic perspective, it is natural on the other hand to classify curve singularities by their genus, notwithstanding the fact that the embedding dimension grows linearly in $b$ - and we have very little control on the deformation theory of singularities of codimension four or more. Isolated Gorenstein singularities of genus one have been classified by Smyth in \cite{SMY1} (with applications to moduli theory in \cite{SMY2,SMY3} and the other references above). Isolated Gorenstein singularities of genus two have been classified by the author in \cite{Bat19}; some non-isolated singularities, that we call ``ribbons with tails'', appear in \cite{BatCar}. Smoothable hyperelliptic singularities have been classified in \cite{BB23}.

A theme has appeared in the author's works on hyperelliptic singularities: namely, the connection between smoothable Gorenstein singularities and logarithmic differentials, i.e. abelian (or more precisely multiscale \cite{BCGGM}) differentials on nodal curves, together with their tropical incarnations:  piecewise-linear functions on the dual graph belonging to the tropical canonical linear series. We outline this connection briefly here, and refer the reader to the next section for more details.
\begin{itemize}[leftmargin=0.5cm]
 \item Given a smoothable Gorenstein curve $C$ of genus $g$ with a unique singular point $p$ of genus $g$, let $\eta_0$ be an abelian differential generating $\omega_{C,p}$. By smoothing $(C,\eta_0)$ and then applying semistable reduction (i.e. the properness of the moduli space of generalised multiscale differentials), we may replace the central fibre with a reducible nodal curve $C^{ss}$ and a multiscale differential $\eta^{ss}$. The preimage $Z$ of $p$ in $C^{ss}$, and the restriction of $\eta^{ss}$ to it, give us a smoothable logarithmic differential.
 \item Given a multiscale differential, we can put it in a one-parameter smoothing, and, possibly after blowing up some points in the central fibre, conjecturally use its various level truncations to produce corresponding Gorenstein contractions. The multiscale differential should then descend to a generator of the dualising bundle of the Gorenstein curves.
\end{itemize}

In line with this program, the main result of this note is the following:
\begin{theorem*}
 There are seven families of isolated Gorenstein singularities of genus three, each corresponding to a connected component of a stratum of abelian differentials of genus three.
\end{theorem*}
Section \ref{sec:sing} consists of the proof of this result; Section \ref{sec:background} describes the geometry of genus three curves, theta characteristics, Gorenstein singularities, and strata of differentials. In the following table we provide a list of the strata and the corresponding singularities, including the gap sequence and minimum number of branches. In the interest of legibility, we have simplified the generators by setting all crimping variables equal to $0$ or $1$ as appropriate. 
We work throughout over $\CC$.
\begin{center}
\begin{tabular}{| c | c | c | c | c | c |}
\hline
 diff. str. & gap seq. & $b_{min}$ & simplified generators & $k:\tm^k\subseteq\m$ & planar\\
 \hline\hline
 
 $\Hcal(1^4)$ & $[2,1]$ & $4$ & \begin{minipage}[t]{6cm} 
  {\small{{$\!\begin{aligned}
  A_3= && \mathbf{c} t_1 && t_2 && t_3 && 0 && 0 && \ldots && 0 \\
  A_4=&& t_1 && \mathbf{c} t_2 && 0 && t_4 && 0 && \ldots && 0 \\
  A_5= && t_1^2 && 0 && 0 && 0  && t_5 && \ldots && 0 \\
  \ldots &&&&&&&&&&&&&& \\
  A_b= && t_1^2 && 0 && 0 && 0  && 0 && \ldots && t_b \\
  B_1= && t_1^2 && t_2^2 && 0 && 0  && 0 && \ldots && 0 \\
 \end{aligned}$} \newline
 w/ analytic invariant $\mathbf{c}\in\CC\setminus\{\pm1\}$}}
  \end{minipage}
 & $\pmod{\tm^3}$ & $4$-uple pt\\ 
 \hline
 
\begin{comment} 
 $\Hcal(1^4)^\text{hyp}$ & $[2,1]$ & $4$ & 
 {\small{{$\!\begin{aligned}
  A_3= && 0 && t_2 && t_3 && 0 && 0 && \ldots && 0 \\
  A_4=&& t_1 && 0 && 0 && t_4 && 0 && \ldots && 0 \\
  A_5= && t_1^2 && 0 && 0 && 0  && t_5 && \ldots && 0 \\
  \ldots &&&&&&&&&&&&&& \\
  A_b= && t_1^2 && 0 && 0 && 0  && 0 && \ldots && t_b \\
  B_2= && t_1^2 && t_2^2 && 0 && 0  && 0 && \ldots && 0 \\
 \end{aligned}$}}}
 & $\pmod{\tm^3}$ & \\ 
 \hline
\end{comment} 
 
 $\Hcal(2,1^2)$ & $[2,0,1]$ & $3$ &  {\small{{$\! \begin{aligned}
   A_3= && 0 && t_2 && t_3 && 0 && \ldots && 0 \\
  A_4= && t_1^3 && 0 && 0 && t_4 && \ldots && 0 \\
  \ldots &&&&&&&&&& \\
  A_b= && t_1^3 && 0 && 0 && 0 && \ldots && t_b \\
  B_1= && t_1^2 && 0 && 0 && 0 && \ldots && 0 \\
  B_2= &&  t_1^3 && t_2^2 && 0 && 0 && \ldots && 0 \\
 \end{aligned}$}}} & $\pmod{\tm^4}$ &\\ 
 \hline
 
 $\Hcal(2^2)^{odd}$ & $[2,0,1]$ & $2$ &  {\small{{$\!\begin{aligned} 
  A_3= && t_1^3 && 0 && t_3  && \ldots && 0 \\
  \ldots &&&&&&&&&& \\
  A_b= && t_1^3 && 0 &&  0 && \ldots && t_b \\
  B_1= && t_1^2 && 0 && 0 &&  \ldots && 0 \\
  B_2= && 0 && t_2^2 && 0 && \ldots && 0 \\
  C_2= && t_1^3 && t_2^3 && 0 && \ldots && 0 \\
 \end{aligned}$}}} & $\pmod{\tm^4}$ & \\ 
 \hline
 
 $\Hcal(2^2)^{ev}$ & $[1,1,1]$ & $2$ &  {\small{{$\!\begin{aligned}
  A_2= && t_1 && t_2 && 0 && \ldots && 0 \\
  A_3= && t_1^3 && 0 && t_3 && \ldots && 0 \\ 
  \ldots &&&&&&&&&& \\
  A_b= && t_1^3 && 0  && 0 && \ldots && t_b \\
  {\color{gray}D_1=} && {\color{gray}t_1^4} && {\color{gray}-t_2^4}  && {\color{gray}0} && \ldots && {\color{gray}0} \\
 \end{aligned}$}}} & $\pmod{\tm^4}$ & $A_7, D_8$\\ 
 \hline
 
 $\Hcal(3,1)$ & $[1,1,0,1]$ & $2$ &  {\small{{$\!\begin{aligned}
  A_2= && t_1^2 && t_2 && 0 && \ldots && 0 \\
  A_3= && t_1^4 && 0 && t_3 && \ldots && 0 \\
   \ldots &&&&&&&&&& \\
  A_b= && t_1^4 && 0 && 0 && \ldots && t_b \\
  C_1= && t_1^3 && 0 && 0 && \ldots && 0 \\
 \end{aligned}$}}} & $\pmod{\tm^5}$ & $E_7$\\ 
 \hline
 
 $\Hcal(4)^{odd}$ & $[1,1,0,0,1]$ & $1$ &  {\small{{$\!\begin{aligned}
 A_2= && t_1^5 && t_2 && 0 && \ldots && 0 \\
   \ldots &&&&&&&&&& \\
 A_b= && t_1^5 && 0 && 0 && \ldots && t_b \\
 C_1= && t_1^3 && 0 && 0 && \ldots && 0 \\
 D_1= && t_1^4 && 0 && 0 && \ldots && 0 \\
 \end{aligned}$}}} & $\pmod{\tm^6}$ & $E_6$\\ 
 \hline
 
 $\Hcal(4)^{ev}$ & $[1,0,1,0,1]$ & $1$ &  {\small{{$\!\begin{aligned}
 A_2= && t_1^5 && t_2 && 0 && \ldots && 0 \\
   \ldots &&&&&&&&&& \\
 A_b= && t_1^5 && 0 && 0 && \ldots && t_b \\
 B_1= && t_1^2 && 0 && 0 && \ldots && 0 \\
 \end{aligned}$}}} & $\pmod{\tm^6}$ & $A_6, D_7$\\ 
 \hline
\end{tabular}
\end{center}

\section{Geometric context}\label{sec:background}

\subsection{The canonical linear series on smooth curves} Let $C$ be a smooth curve of genus three. The canonical linear series is a basepoint-free $\mathfrak{g}^2_4$ on $C$, producing it either as a smooth plane quartic, or as a double cover of the Veronese conic (if and only if $C$ is hyperelliptic).

For a general point $P$ of $C$, the orders of vanishing of $\lvert K_C\rvert$ at $P$ are $0<1<2$. When the vanishing sequence is different from this, i.e. $P$ is a ramification point of the canonical linear series, $P$ is called a \emph{Weierstrass point}. If the vanishing orders are denoted by $a_0<a_1<a_2$, the orders of ramification form the non-decreasing sequence $(\alpha_i=a_i-i)_{i=0,1,2}$. The weight of a Weierstrass point is defined as $\sum\alpha_i$. The classical Pl\"ucker formula states that the sum of the weights of all Weierstrass points is $(g-1)g(g+1)=24$ (see e.g. \cite[Chapter I, Exercise C]{ACGH}).

If $C$ is hyperelliptic, the ramification sequence of a Weierstrass point is $(0,1,2)$, hence there are exactly $8$ Weierstrass points.\footnote{This can also be deduced from the Riemann--Hurwitz formula.} Ignoring automorphisms, the hyperelliptic locus can be thus identified with $\Mcal_{0,8}/\mathfrak{S}_8$. It forms a divisor in $\Mcal_3$ whose class was computed by Mumford at the very beginning of the subject of intersection theory on moduli spaces \cite{MumfordEnumGeom} (the class of its closure in $\oM_3$ is also well-known, see \cite[Chapter 3, \S H]{HM} and \cite{EstevesPorteous}).

If $C$ is a general plane quartic, it has exactly $24$ flexes, i.e. simple Weierstrass points with ramification sequence $(0,0,1)$. Interestingly, the general plane quartic can be reconstructed from its inflection lines (even without knowing where they meet the curve!), see \cite{PaciniTesta}. A hyperflex is a point of ramification sequence $(0,0,2)$, hence there are at most $12$ of them. The locus of curves with a hyperflex is an irreducible divisor in the space of quartic plane curves, whose class was studied by Faber \cite{FaberChow}. In fact, a smooth quartic with any number $n$ of hyperflexes, $0\leq n\leq9$ or $n=12$, can be found; see Vermeulen's thesis \cite{Vermeulen}. 
A stratification of plane quartics according to their singularity type can be found in \cite{HuiThesis}.

\subsection{Theta characteristics and singular curves} A theta characteristic $L$ is a  square-root of the canonical line bundle. The space of theta characteristics $S\subseteq\on{Pic}^{g-1}(C)$ is a torsor under the $2$-torsion subgroup of the Jacobian $J_2(C)$; in particular, if $C$ is smooth, $S$ consists of $2^{2g}$ elements. It is a classical fact due to Riemann, Atiyah \cite{Atiyah}, and Mumford \cite{MumfordTheta} that the \emph{parity} of the dimension $h^0(L)$ of the space of regular sections of a theta characteristic is locally constant in families. The space $S$ of theta characteristics can be divided into two subspaces $S^+$ and $S^-$ accordingly. For a fixed smooth curve of genus $g$,
\[ \# S^+=2^{g-1}(2^g+1),\qquad \# S^-=2^{g-1}(2^g-1),\]
summing up to $2^{2g}$.

\subsubsection{Bitangents} Bitangents of a smooth quartic curve are odd theta characteristics, so there are $28$ of them. In general, the $36$ even theta characteristics are not effective. In the hyperelliptic case $\psi\colon C\to \PP^1$, instead, there is an effective even theta characteristic corresponding to the class of $\psi^*\OO_{\PP^1}(1)$; the $28$ odd theta characteristics can be written as $q_1+q_2$, with $q_1$ and $q_2$ two distinct Weierstrass points.
Any smooth quartic can be reconstructed from its bitangent lines; in fact from any of a number of special configurations of seven of them, known as Aronhold systems \cite{CaporasoSernesi,Lehavi}. 

\subsubsection{Projective duality} Flexes and bitangents play an important role in projective duality. Let $C$ be a reduced plane curve of degree $d$. The \emph{projectively dual} curve $C^*$  is the closure of the locus of tangent lines to smooth points of $C$ in the dual projective plane $\PP^{2*}$. If $C$ is smooth, $C^*$ has degree $d^*=d(d-1)$, but it is usually singular: indeed, a bitangent of $C$ corresponds to a node of $C^*$, and an ordinary flex to a cusp. Let $\delta^*$ (resp. $\kappa^*$) denote the number of bitangents (resp. flexes) of $C$. If $C$ is general, the Pl\"ucker formulae imply:
\begin{align*}
 \delta^*=&\frac{1}{2}(d+3)d(d-2)(d-3),\\
 \kappa^*=& 3d(d-2).
\end{align*}
For quartics we recover: $28$ bitangents and $24$ flexes.
When $C$ is singular, the formulae above need to be corrected by taking the singularities into account. See the work of Wall \cite{Wall} on how singularities behave under projective duality; for instance, a hyperflex corresponds to an $E_6$-singularity (of genus three, see below).

\subsubsection{Harris}\label{sec:Harris} Theta characteristics can be defined as above on any Gorenstein curve \cite{HarrisTheta} (more generally, one needs to consider rank-one torsion-free sheaves, because such is the dualising sheaf itself \cite{Piontkowski}). The parity of the space of global sections remains locally constant, and the space of theta characteristics remains a torsor under the $2$-torsion subgroup of the Jacobian $J_2(C)$, but the cardinality of the latter may vary. Indeed, let $C^{\nu}\to C^{sn}\to C$ denote the normalisation (respectively, seminormalisation, i.e. the initial curve homeomorphic to $C$) of the reduced curve $C$. The kernel of $\Pic(C^{sn})\to\Pic(C^{\nu})$ is an algebraic torus of dimension $\beta$, %=b_1(\Gamma)$, the first Betti numer of the dual graph of $C$\footnote{The dual graph can be defined in the same way as for nodal curves, i.e. it has one vertex per irreducible component of $C$, with one edge connecting two vertices for every point of intersection of the corresponding components, independent of its local structure.};
so the cardinality of its $2$-torsion is $2^\beta$. On the other hand, the kernel of $\Pic(C)\to\Pic(C^{sn})$ is unipotent, i.e. an extension of copies of $\Ga$; in particular it is divisible and it has trivial $2$-torsion. So, $J_2(C)$ has cardinality $2^{2\tilde g+\beta}$, where $\tilde g$ is the sum of the genera of the components of $C^\nu$. See for instance \cite[\S 7.5]{Liu}.

We offer a brief summary of (Mumford's and) Harris' work: given a theta characteristic $L$, the function
\[q_L(D)=h^0(L(D))-h^0(L) \pmod 2\]
is a quadratic form on the $\ZZ/2\ZZ$-vector space $J_2(C)$; its associated bilinear form is the Weil pairing.

Let $\Gamma_2$ denote the kernel of the pullback map $J_2(C)\to J_2(C^\nu)$. The restriction of $q_L$ to $\Gamma_2$ is linear and independent of $L$; denote it by $l$. Interestingly, $l$ is determined by local data. Indeed, we can describe a set of generators of $\Gamma_2$ as follows: we let $e_i$ be the line bundle that is trivial on $C^\nu$, and glued with $-1$ between branch $i$ and any other branch, and $1$ between any two branches other than the $i$-th. These line bundles descend to $C$, and they satisfy exactly one relation $\sum e_i=0$. If $D$ denotes an adjoint divisor, that is defined by the pullback of the conductor ideal to $C^\nu$, Harris proves that \[l(e_i)=\on{mult}_{p_i}(D),\] where $p_i$ are the points of $C^\nu$ over $p\in C$.

If $l$ is not $0$, then the set $S$ of theta characteristics subdivides equally between even and odd thetas. If $l\equiv 0$, instead, $q_L$ is induced by a quadratic form on $J_2(C^\nu)$, and there are two possibilities for the latter: in particular, the number of even - resp. odd - theta characteristics is $2^{\tilde g+\beta-1}(2^{\tilde g}+1)$ - resp. $2^{\tilde g+\beta-1}(2^{\tilde g}-1)$ -, or viceversa. Let us set $Q(C)=0$ in the former case, and $Q(C)=1$ in the latter. Once again, Harris proves that this dichotomy is determined only by local data around the singularities of $C$. Indeed, if $l\equiv0$, the adjoint divisor $D$ is even, i.e. $D=2E$; let $I'$ denote the ideal of $\OO_C$ that pulls back to $\OO(-E)$ on $C^\nu$, and let $\epsilon(p)=\on{len}(\OO_C/I')$. Then \[Q(C)=\sum_{p\in C^{\text{sing}}} \epsilon(p) \pmod 2.\]

Notice that the space of theta characteristics is at any rate finite (albeit not \`etale) over the moduli stack of Gorenstein curves.

\subsection{Isolated Gorenstein curve singularities}\label{sec:Gorenstein}

A Cohen--Macaulay variety is called Gorenstein if its dualising sheaf is a line bundle. This is a local condition that can be checked after taking completions. Here we shall study analytic germs of curve singularities. For example, plane curve (and more generally l.c.i.) singularities are Gorenstein. Gorenstein curves need not be reduced: being Cohen--Macaulay, they have no embedded points, but they may still contain multiple components; in this case, the generic slice\footnote{Locally at the generic point of any component, the underlying reduced curve is regular and defined by one parameter. By quotienting this out, we obtain a $0$-dimensional Gorenstein ring.} must itself be Gorenstein (for example, double structures always work, but, among the triple structures, only the curvilinear ones are allowed). However, we will focus on isolated singularities in this paper. 

\subsubsection{Rosenlicht}\label{par:Rosenlicht}
Recall that the dualising bundle of a smooth variety is identified with the determinant of its cotangent bundle, also called the canonical bundle. Rosenlicht gave an explicit description of the dualising sheaf of a reduced curve $C$ in terms of meromorphic differentials on its normalisation $C^\nu$ (see for instance \cite[Chapter VIII]{AK}): a rational differential (on the normalisation) $\eta \in \omega_{K(C)/\CC}$ descends to a regular section $\eta\in \Gamma(C,\omega_C)$ if and only if, for every (germ of) regular function $f\in \OO_{C,p}$, the sum of the \emph{residues} of $f\cdot\eta$ at the points of $C^{\nu}$ over $p$ is zero:
\[\sum_{p_i\in\nu^{-1}(p)}\on{Res}_{p_i}(f\cdot\eta)=0.\]

Let $\cond=\on{Ann}(\nu_*\OO_{C^\nu}/\OO_C)$ denote the conductor ideal of the singularity; it can be viewed both as an ideal of $\OO_C$ and as one of $\nu_*\OO_{C^\nu}$. If $D$ is the associated (adjoint) divisor on $C^\nu$, notice that $\on{mult}_{p_i}(D)$ bounds from above the possible orders of pole of $\eta$ at $p_i$; in particular, $\eta$ must be regular at smooth points of $C$. Indeed, the residue pairing:
\begin{equation}\label{eq:residue}
(\nu_*\OO_{C^\nu}/\OO_C)\times(\omega_C/\nu_*\omega_{C^\nu})\to\CC
\end{equation}
is perfect, and when looking at germs it can be enhanced to a graded pairing.

Let $\delta=\on{len}(\nu_*\OO_{C^\nu,p}/\OO_{C,p})$; this is called the \emph{$\delta$-invariant} of $(C,p)$. It represents the number of conditions that a function on the normalisation must satisfy in order to descend to the singularity. If we factor the normalisation $\nu$ through the seminormalisation, $C^\nu\to C^{sn}\to C$, and if $b$ is the number of branches of $C$ at $p$ (i.e. the number of analytically connected components of the normalisation), then there are exactly $b-1$ conditions for a function to descend to $C^{sn}$, expressing the requirement that functions attain the same value at the $b$ preimages of $p$ in $C^\nu$. Indeed, the germ of the seminormalisation is isomorphic to the union of the coordinate axes in $\Aaff^b$. Dually under \eqref{eq:residue} we find the differentials $\frac{\de t_1}{t_1}-\frac{\de t_i}{t_i},\ i=2,\ldots,b,$ on the seminormalisation; in particular, the seminormalisation is never Gorenstein for $b>2$.

The remaining \[g=\delta-b+1\]
conditions determine the so-called \emph{genus} of the singularity. This is also the local contribution of the singularity to the arithmetic genus of a projective curve containing it. For instance, the cusp and the tacnode are singularities of genus $1$; they appear naturally in families of plane cubics.

\subsubsection{Serre}\label{sec:Serre} For graded finite-dimensional $\CC$-algebras, the Gorenstein condition reduces to Poincar\'e duality. More generally, a ring is Gorenstein if it can be cut down to a Gorenstein Artinian ring by means of a regular sequence.

A monomial unibranch curve is Gorenstein if and only if its gap sequence is symmetric: if $J=\{j_1,\ldots,j_g\}$ are the exponents of the missing monomials, then $n\in J$ if and only if $2g-1-n\notin J$ \cite{Kunz}.

For reduced curves, Serre \cite[\S 4, Proposition 7]{Ser88} shows that an equivalent condition to being Gorenstein is
\begin{equation}\label{eqn:Serre}
 \on{len}(\OO_{C,p}/\cond)=\on{len}(\nu_*\OO_{C^\nu,p}/\OO_{C,p}) (=\delta).
\end{equation}

\subsubsection{Stevens}\label{sec:Stevens} A curve singularity $(C,p)$ is \emph{decomposable} if $C$ is the union of two curves $C_1,\ C_2$, living in smooth spaces that intersect transversely in $p$. Equivalently, there are two parameters $t_1$ and $t_2$ on the normalisation, and a partition $(x_1,\ldots,x_a,y_1,\ldots,y_b)$ of the generators of $\m/\m^2$, such that the $\mathbf x$ do not involve $t_2$, and the $\mathbf y$ do not involve $t_1$. A useful observation \cite[Proposition 2.1]{AFSGm} is that every Gorenstein curve singularity, except for the node, is indecomposable. For instance, the rational (or ordinary) $b$-fold point, i.e. the singularity of the seminormalisation, is not Gorenstein for $b\geq3$.

\subsubsection{Smyth}\label{sec:Smyth} Let $(R,\m)$ denote the complete local ring of an isolated curve singularity, with normalisation $(\tR,\tm)\simeq\left(\CC[\![t_1]\!]\oplus\ldots\oplus\CC[\![t_b]\!],\langle t_1,\ldots,t_b\rangle\right)$. Following \cite[Appendix A]{SMY1} and \cite[\S 1]{Bat19}, we consider $\tR/R$ as a $\mathbb Z$-graded $R$-module (and finite $\CC$-algebra) with:
\[ (\tR/R)_i:=\tm^i/(\tm^i\cap R)+\tm^{i+1}.\]
The following observations provide the starting point for the classification:

\begin{enumerate}
\item $b+2=\delta(p)=\sum_{i\geq 0}\dim_\CC(\tR/R)_i;$
\item $3=g=\sum_{i\geq 1}\dim_\CC(\tR/R)_i;$
\item\label{obs:add} if $(\tR/R)_i=(\tR/R)_j=0$ then $(\tR/R)_{i+j}=0$;
\item\label{obs:igrad} $\sum_{i\geq j}(\tR/R)_i$ is a grading of $\tm^j/(\tm^j\cap R)$;
\item\label{obs:ses} there are short exact sequences of $R/\m=\CC$-modules:
 \[ 0\to K_i:=\frac{\tm^i\cap R}{\tm^{i+1}\cap R}\to \frac{\tm^i}{\tm^{i+1}}\to \left(\tR/R\right)_i\to 0.\]
\end{enumerate}

Accordingly, the set \[\{i\in\mathbb{N}\ |\ (\tR/R)_i=0\}\cup\{0\}\] is a submonoid; to the remaining integers we may associate the dimension of $(\tR/R)_i$: this forms a non-increasing sequence that we call the \emph{gap sequence}.

\subsection{Semistable tails and differentials}\label{sec:sstails}
Suppose that a Gorenstein singularity $(C,p)$ is smoothable. Consider it as a projective curve $C$ by compactifying every branch to a distinct, smooth and disjoint away from $p$, copy of $\PP^1$. Then $C$ is smoothable as well by \cite[Example 29.10.2]{HarDef}. Consider a one-parameter smoothing $\mathcal C\to\dvr$, and apply semistable reduction to $\mathcal C$, in such a way to obtain a (minimal) contraction $\Phi\colon\mathcal C^{ss}\to\mathcal C$. We now focus on the special fibre. What we see is a nodal curve $C^{ss}$ of genus $g$ with a morphism $\phi$ to $C$. Let $Z$ denote the preimage of $p$. Then $C_0=\overline{C\setminus Z}$ is the disjoint union of $b$ copies of $\PP^1$, and $\phi$ presents it as the normalisation of $C$. The exceptional curve $Z$, instead, is a nodal curve of genus $g$, naturally marked by its intersection with the $b$ branches of $C_0$. It is called a \emph{semistable tail} of the singularity. Although $Z$ does depend on the chosen smoothing, we shall see that we can provide a uniform characterisation in terms of logarithmic differentials.

Let $\eta$ be a local generator of the dualising sheaf $\omega_{C,p}$. Since the singularity is indecomposable, the multiplicity of the conductor ideal at every preimage $p_i$ in the normalisation is at least two, and so is the order of pole of $\eta$. This means that, when compactifying the $i$-th branch to $\PP^1$ via $s_i=1/t_i$, the (meromorphic) differential $\eta_i=\on{d}\!t_i/t_i^{m_i}+h.o.t.$ ($\eta$ restricted to the $i$-th branch $C_i$) extends to a holomorphic differential on the complement of $p_i$. If $m_i$ is the order of pole at $p_i$, it will indeed have $m_i-2$ zeroes on this branch. Thus $\eta\in\Gamma(C,\omega_C)$, and $\on{deg}(\omega_{C|C_i})=m_i-2$.

Going back to the smoothing family, consider the relative dualising bundle $\omega_{\Ccal/\dvr}$. Its pullback $\Phi^*\omega_{\Ccal/\dvr}$ to $\Ccal^{ss}$:
\begin{itemize}
 \item is trivial on the exceptional locus $Z$;
 \item coincides with $\omega_{\Ccal^{ss}/\dvr}$ outside $Z$, where $\Phi$ is an isomorphism.
\end{itemize}
Assuming regularity of $\Ccal^{ss}$, then, we may write $\Phi^*\omega_{\Ccal/\dvr}=\omega_{\Ccal^{ss}/\dvr}(\widetilde{Z})$, where $\widetilde{Z}$ is a (vertical) Cartier divisor supported on $Z$. On the other hand, the previous paragraph implies that the multiplicity of $\widetilde{Z}$ at the $i$-th component of $C_0$ is $m_i-1$ (the $1$ comes from Noether's formula:  for a nodal curve $X$, and a closed subcurve $Y\subset X$, the restriction of $\omega_Y$ to $X$ is $\omega_{X|Y}=\omega_Y(Y\cap\overline{X\setminus Y})$; more generally, the difference is an adjoint divisor, see \cite[Proposition 1.2]{Catanese}). We conclude: \[\omega_{Z}\left(\sum_{i=1}^b (2-m_i)p_i\right)=\OO_{Z},\]
that is $\sum_{i=1}^b (m_i-2)p_i$ is a canonical divisor on $Z$.\footnote{If $Z$ is reducible, these multiplicities determine further the shape of the dual graph: for instance, when $g=1$, all the $p_i$ must be at the same distance from the minimal genus one subcurve \cite{SMY1}. In general, these conditions can be phrased effectively in terms of logarithmic geometry, by requiring the existence of certain piecewise-affine functions on the dual graph, which turn out to be differentials in the sense of tropical geometry: see \cite{Bat19} and \S\ref{sec:differentials} below. It was pointed out to me by F. Viviani that this suggests the following statement: the variety of stable tails \cite{HassettStableTails} of a Gorenstein singularity is the Deligne--Mumford compactification of a stratum of differentials.}

The multiplicities $m_i-2$ are determined by the singularity $(C,p)$ via its conductor. This determines a stratum of holomorphic differentials. Based on the classification of genus three singularities, we conjecture that every topological type of Gorenstein singularity determines in fact \emph{a connected component} of a stratum of holomorphic differentials. We wonder whether the converse is also true.

\subsection{(Compact) strata of holomorphic differentials}\label{sec:differentials}
The moduli space of holomorphic differentials $\Hcal$ on smooth projective curves is stratified by the multiplicities of the zeroes $\mu=(m_i)\vdash 2g-2$. Every stratum consists of at most three connected components, depending on $\mu$ \cite{KontsevichZorich}:
\begin{itemize}[leftmargin=.5cm]
 \item In case $\mu=(2g-2)$ or $\mu=(g-1,g-1)$, there is a hyperelliptic stratum where the underlying curve is hyperelliptic, the differential is anti-invariant with respect to the hyperelliptic involution,\footnote{A necessary condition for an anti-invariant differential with prescribed multiplicity $\mu$ to exist is that every odd part $m_i$ appears an even number of times. Except for the two cases listed above, these hyperelliptic differentials arise as limits of differentials on smooth curves.} and the zero(es) are located at a Weierstrass (resp. two conjugate) points.
 \item If every part $m_i$ of $\mu$ is even, then $\Theta=\sum\frac{m_i}{2}p_i$ is an effective theta characteristic, and there are two connected components according to the parity of $h^0(\Theta)$ (the underlying curve is general).
\end{itemize}
 
\begin{rem}
 In genus three, the presence of an effective even theta characteristic implies that the underlying curve is hyperelliptic. Hence there are always at most two connected component for every stratum when $g=3$.
\end{rem}

There are several ways of compactifying strata of holomorphic differentials. The underlying curve $C$ is usually allowed to become nodal and reducible. If we consider differentials up to scaling, then:
\begin{itemize}[leftmargin=.5cm]
 \item they are sections of the projectivisation of the Hodge bundle $\mathbb{E}=\pi_*(\omega_{\Ccal/\Mcal_g})$;
 \item they are determined by their vanishing locus, i.e. the associated divisor $\sum m_ip_i$.
\end{itemize}
Correspondingly, we can compactify $\Hcal(\mu)$ by taking its closure in $\PP(\mathbb E)$ - notice that the Hodge bundle extends naturally to the Deligne--Mumford compactification -, resp. in $\oM_{g,b}$, where it has codimension $g-1$. Unfortunately, taking closures is problematic from the point of view of moduli theory.
The moduli space of multiscale differentials provides instead a modular compactification  \cite{BCGGM}, admitting a log geometric interpretation \cite{Chen2,Tale}; here is a naive description.

Consider a semistable one-parameter family $\Ccal\to\dvr$ as above, and suppose the general fibre $C_\xi$ is smooth, and $\eta_\xi$ is a holomorphic differential on it. The central fibre $C_0$ may be nodal and reducible. The differential $\eta_\xi$ can be extended to a differential $\eta$ on the whole family, but the latter may vanish on a subcurve $C_0^\prime$ of $C_0$ (equivalently, the associated divisor may have a non-trivial vertical part). Let $C_{0,0}$ be the subcurve on which $\eta$ does not vanish, and let $\eta_0$ denote its restriction. Rescaling the differential by an appropriate (negative) power of the base parameter $t$, though, it is possible to extend $\eta$ to a non-trivial differential on $C_0^\prime$; this will have poles along $C_{0,0}$. Let $\eta_1$ denote the meromorphic restriction of $\eta^\prime$ to $C_0^\prime$, and $C_{0,-1}\subseteq C_0^\prime$ the subcurve on which it is not trivial, and so on.

We end up with a collection $\{\eta_i\}$ of meromorphic differentials on subcurves $\{C_{0i}\}$ partitioning the central fibre. Moreover, the dual graph of the latter appears now to be ordered according to the order of vanishing of the differential $\eta$: we obtain a so-called level graph by putting $C_{0i}$ at level $i$. This information can be enhanced to a compatible integral piecewise-affine function $\lambda$ on the dual graph, the slopes recording the order of vanishing/poles of the differential $\eta_i$ at the corresponding special point of the subcurve $C_{0i}$ (as we have seen above, Noether's formula implies that the slope differs by the actual order by $\pm1$). In the logarithmic interpretation of multiscale differentials, $\lambda$ is nothing but the tropicalisation of $\eta$. All of this is, roughly speaking, what is called a \emph{generalised multiscale differential}.

Unfortunately, the space of generalised multiscale differentials is typically not a topological compactification. In fact, it is not irreducible, but the locus of smoothable differentials has been characterized in terms of a \emph{global residue condition} \cite{BCGGM16}: a zero-sum condition on residues of $\eta_i$ at poles belonging to possibly different irreducible components of the curve, which are however connected through components at higher levels. Here we reproduce a conjecture - originally due to D. Ranganathan and J. Wise, and spelled out in \cite{BB23} - relating the smoothability of multiscale differentials and Gorenstein curves in purely algebraic (and characteristic-free) terms.

It turns out that the piecewise affine function $\lambda$ above (or rather its negative) is indeed an element of the tropical canonical linear series. The same is true for every truncation of $\lambda$ to the subcurve of $C_0$ lying at level $i$ or below. Even for these purely combinatorial data, the moduli space is not in general irreducible (nor pure-dimensional); the locus of smoothable (\emph{realisable}) tropical differentials has been described explicitly in \cite{MUW}. We thus have the first part of:

\begin{conj}[differential descent]\label{conj} \leavevmode
Let $(C,\eta)$ be a generalised multiscale differential with tropicalisation $\lambda$. Then $\eta$ is smoothable if and only if \begin{enumerate}[label=(\roman*)]
    \item for every level $i$, the truncation $\lambda_i$ of $\lambda$ is a \emph{realisable} tropical differential;
    \item there exists a logarithmic modification $\widetilde{C}\to C$, a natural extension $\tilde\eta$ of the pullback of $\eta$ to $\widetilde C$, and a reduced \emph{Gorenstein} contraction $\sigma\colon \widetilde C\to\oC_i$ such that $\sigma^*\omega_{\oC_i}=\omega_C(\lambda_i)$, and
    \item the differential $\tilde\eta_i$ at level $i$ descends to a local \emph{generator} of $\omega_{\oC_i}$.
\end{enumerate}
\end{conj}

We refer the reader to \cite{BB23} for a more detailed discussion, particularly of the hyperelliptic case.

\subsection{Crimping spaces and moduli of singular curves}\label{sec:crimping}
We think of a curve singularity as being obtained from the normalisation by pinching a finite subscheme \cite[Lemma 1.2]{vdW}. Algebraically, this presents the local ring as a fiber product, contained in the local ring of the normalisation. The \emph{crimping space} (or moduli of attaching data) can be thought of as the moduli space of finite subschemes of $C^\nu$ whose pinching produces a singularity with fixed analytic germ (or type $\tau$), see for instance \cite[\S 2.2]{SMY2}. Algebraically, it is the moduli space of subalgebras $R$ of type $\tau$ in $\tR=\CC[\![t_1]\!]\oplus\ldots\oplus\CC[\![t_b]\!]$, see \cite[\S 1.3]{vdW}. Since the conductor $\cond$ is contained in $R$, and the numerics of the type determine a power $k$ of $\tm$ that is contained in $\cond$, $R$ is completely determined by its image in the Artinian quotient $\tR/\tm^k$. In this way we may view the crimping space as the quotient of $\on{Aut}_{\CC}(\tR/\tm^k)$ by the stabiliser of $R/\tm^k$, and in particular, by a result of Rosenlicht, it is an extension of an affine space $\Aaff^h$ with a product of $\Aaff^1\setminus\{0\}$ \cite[Theorem 1.74]{vdW}.

Here is some heuristic motivation of why crimping spaces are important. There are several examples of birational models of $\oM_{g,n}$ allowing for more singular curves in exchange for nodal subcurves with few special points, see the introduction for a sample. Crimping spaces and their compactifications appear as fibres in the flipping loci. From the above description we gather that they are rational, and that they can be described in terms of differential-geometric information on the normalisation; in particular, their classes are tautological. This is particularly relevant to the enumerative applications.

\section{Algebraic classification}\label{sec:sing}
We organise the classification according to the finitely many possible gap sequences, see \S \ref{sec:Gorenstein}. If $i$ is a gap of dimension $a_i$, we write $[a_i]_{i\geq1}$ for the vanising sequence, stopping when $a_i$ becomes definitely $0$. We look for generators of $\m/\m^2$. We write $A_i=\sum A_i(j,k)t_j^k$ for the linear generators (there exists a $j$ such that $A(j,1)\neq0$), $B_i$ for the quadratic ones, etc.
If $i$ is the last non-trivial gap, then $\tm^{i+1}\subseteq R$, so in particular $\tm^{i+1}\subseteq \cond$. This means that we only have to look for generators modulo $\tm^{i+1}$. We adopt two strategies simultaneously: studying the dimension of $R/\cond$ as a $\CC$-vector space, cf. Formula \eqref{eqn:Serre}; and writing down the relations that a generator $\eta$ of the dualising sheaf must satisfy, cf. Formula \eqref{eq:residue}.

\subsection*{[3]} Cannot possibly be Gorenstein, since $\tm^2\subseteq\frakc$, so:
 \[\dim_{\CC}(R/\frakc)\leq\dim_{\CC}(R/\tm^2)=\dim_{\CC}(K_0+K_1)=1+(b-3)<\delta.\]
 
\subsection*{[2 1]} Note that $\tm^3\subseteq R$. The Gorenstein condition implies $b\geq 4$: indeed, if $b=3$, $R/\cond$ would be generated by $1$, at most one linear generator $A$, and two quadratic ones $B_1$ and $B_2$, hence $5=\delta=\dim R/\cond\leq 4$, a contradiction (similarly for $b=1,2$).
 
 Up to relabelling, assume that $t_1^2$ does not belong to $R$. Then we have quadratic elements \[B_i=t_i^2+B_i(1,2)t_1^2\] in $R$ for $i=2,\ldots,b$. Then $t_1$ does not belong to $R$ either; there is another linear monomial missing, which, up to relabelling, we may assume is $t_2$. Then we have linear generators \[A_i=t_i+A_i(1,1)t_1+A_i(1,2)t_1^2+A_i(2,1)t_2\] for $i=3,\ldots,b$ (we may reduce to this form by summing appropriate scalar combinations of the $B_i$). By comparing quadratic terms in the $A_i$ with the $B_j$, and since we assumed $t_1^2\notin R$, we find the following relations must hold:
 \begin{subequations}
 \begin{align}
  A_i^2=B_i+A_i(2,1)^2B_2 &\Leftrightarrow & A_i(1,1)^2=B_i(1,2)+A_i(2,1)^2B_2(1,2),\label{eq1}\\
  A_iA_j=A_i(2,1)A_j(2,1)B_2 &\Leftrightarrow & A_i(1,1)A_j(1,1)=A_i(2,1)A_j(2,1)B_2(1,2).\label{eq2}
 \end{align}  
 \end{subequations}
 
 Suppose that $B_2(1,2)=0$, or that $A_i(2,1)=0$ for all $i$. Then, by \eqref{eq2}, there exists at most one index (say $i=3$) such that $A_i(1,1)\neq0$. Then, $B_j=A_j^2=t_j^2$ for $j=4,\ldots,b$, hence $R/\cond\subseteq\langle1,A_3,\ldots,A_b,B_2,A_3^2,B_3\rangle$, with the last three linearly dependent by \eqref{eq1}, so $R/\cond$ has dimension $\leq b+1<\delta$, a contradiction. So we may assume $A_3(2,1)\neq0$ and $B_2(1,2)\neq0$. In particular, $t_2^2$ does not belong to $R$.
 
 Now suppose that $A_j(1,1)=0$ for all $j=4,\ldots,b$. Then, by setting $i=3$ in \eqref{eq2}, we see that $A_j(2,1)=0$ for all $j=4,\ldots,b$ as well, so $B_j=A_j^2=t_j^2$ which gives a contradiction as above. So assume $A_4(1,1)\neq0$.
 
 From \eqref{eq2} we find
 \begin{align*}
  A_4(2,1)&=\frac{A_3(1,1)}{B_2(1,2)A_3(2,1)}A_4(1,1), & \text{and}\\
  A_j(2,1)&=\frac{A_3(1,1)}{B_2(1,2)A_3(2,1)}A_j(1,1) &=\frac{A_3(1,1)A_4(2,1)}{A_3(2,1)A_4(1,1)}A_j(2,1)\\
  &=\frac{A_3(1,1)^2}{A_3(2,1)^2B_2(1,2)}A_j(2,1), &j=5,\ldots,b.
 \end{align*}
 Let
 \begin{equation}\label{crossratio} 
   \Delta:=\frac{A_3(1,1)A_4(2,1)}{A_3(2,1)A_4(1,1)}=\frac{A_3(1,1)^2}{A_3(2,1)^2B_2(1,2)}.
  \end{equation}
 If $\Delta=1$, in particular $A_3(1,1)\neq0$, then $A_j(2,1)=\frac{A_3(2,1)}{A_3(1,1)}A_j(1,1)$, and \[A_j^2-\frac{A_j(1,1)^2}{B_2(1,2)}B_2=t_j^2\] belongs to $R$, therefore to $\cond$, for all $j=4,\ldots,b$. This is a contradiction as above.
 
 So $\Delta\neq 1$, which implies that $A_j(1,1)=A_j(2,1)=0$ for all $j=5,\ldots,b$.
Finally, $A_j(1,2)\neq0$ for $j=5,\ldots,b$, otherwise the singularity would be decomposable.

\noindent Summing up, modulo $\tm^3$, the generators of $\m/\m^2$ can be written as:
 \begin{align*}
  A_3= && A_3(1,1)t_1+A_3(1,2)t_1^2 && A_3(2,1)t_2 && t_3 && 0 && 0 && \ldots && 0 \\
  A_4=&& A_4(1,1)t_1+A_4(1,2)t_1^2 && A_4(2,1)t_2 && 0 && t_4 && 0 && \ldots && 0 \\
  A_5= && A_5(1,2)t_1^2 && 0 && 0 && 0  && t_5 && \ldots && 0 \\
  \ldots &&&&&&&&&&&&&& \\
  A_b= && A_b(1,2)t_1^2 && 0 && 0 && 0  && 0 && \ldots && t_b \\
  B_2= && B_2(1,2)t_1^2 && t_2^2 && 0 && 0  && 0 && \ldots && 0 \\
 \end{align*}
 with $A_3(2,1), A_4(1,1), A_5(1,2),\ldots,A_b(1,2), B_2(1,2)\in\Gm$ and $A_3(1,1), A_4(2,1)$, $A_3(1,2), A_4(1,2)\in\Ga$, such that \[A_4(2,1)=\frac{A_4(1,1)}{A_3(2,1)B_2(1,2)}A_3(1,1),\qquad\text{and}\qquad\Delta=\frac{A_3(1,1)A_4(2,1)}{A_3(2,1)A_4(1,1)}\neq 1.\]
 $R/\cond$ is spanned by $1,A_3,\ldots,A_b$,  $B_2,A_3^2, A_4^2$, with dimension $b+2=\delta$.
  
  The dualising line bundle is generated by the meromorphic differential:
\[\frac{\de t_1}{t_1^3}-B_2(1,2)\frac{\de t_2}{t_2^3}
+A_3(2,1)^2B_2(1,2)(1-\Delta)\frac{\de t_3}{t_3^3}+A_4(1,1)^2(\Delta-1)\frac{\de t_4}{t_4^3}-\sum_{i=3}^m A_i(1,2)\frac{\de t_i}{t_i^2}.
\]

It follows that, if $p_i$ denotes the attaching point of the branch with coordinate $t_i$ to the semistable tail $Z$ (cf. \S\ref{sec:sstails}), then we can write:
\[\omega_Z=\OO_Z(p_1+p_2+p_3+p_4).\]
Thus, the pointed semistable tail $(Z;p_1,\ldots,p_4)$ belongs to the generic (open) stratum of differentials $\Hcal(1,1,1,1)$.
  
\begin{remark}
 In case $A_3(1,1)\neq 0$, the element $B_2$ is a non-zero multiple of $A_3A_4$. For $b=4$, set $x=A_3,\ y=A_4$, and let $a_1,a_2$ and $b_1,b_2$ denote their linear coefficients. The equation of this singularity is $xy(b_1x-a_1y)(b_2x-a_2y)$. The conditions that all the coefficients are non-zero and $\Delta\neq1$ ensure that the four lines are distinct, so that the resulting singularity is reduced. The planar quadruple point is the singularity of smallest genus and number of branches that has analytic moduli, in the sense that the cross-ratio of the four lines is invariant under analytic isomorphisms. See \cite[Example 3.43.2]{Greuel}. This family is a generalisation of the planar quadruple point to an arbitrary number of branches/embedding dimension.
\end{remark}

\begin{remark}
  Once the invariant \eqref{crossratio}, i.e. the analytic type of the singularity, is fixed, the remaining coefficients are parametrised by $\Gm^{b-1}\times\Ga^2$ (as a variety, not as a group). This is the crimping space (see \S \ref{sec:crimping}). It can be described more geometrically in terms of pinching data as follows (set $b=4$ for simplicity). Start from the rational $4$-fold point. Firstly, collapse a generic $2$-plane in the tangent space at the origin, i.e. one that does not intersect any of the coordinate $2$-planes. We can represent it as a map:
 \[\CC[\![t_1,\ldots,t_4]\!]/(t_it_j)\to\CC[\epsilon_1,\epsilon_2]/(\epsilon_1,\epsilon_2)^2.\] By genericity, and up to the $GL_2$-action, the map is represented by the matrix:
 \begin{equation*}
 \begin{pmatrix}
 1 & 0 & -a_1 & -a_2\\
 0 & 1 & -b_1 & -b_2
 \end{pmatrix}
 \end{equation*}
Algebraically, pinching corresponds to the fibre product, with generators:
\begin{align*}
  x^\prime= && a_1t_1 && a_2t_2 &&  t_3 && 0  \\
  y^\prime= && b_1t_1 && b_2t_2 && 0 && t_4  \\
  z^\prime= && t_1^2 && 0 && 0 && 0 \\
 \end{align*}
Secondly, collapse a tangent line not contained in the $(x^\prime,y^\prime)$-plane; we can write the map $\CC[\![x^\prime,y^\prime,z^\prime]\!]/(\ldots)\to\CC[\epsilon^\prime]/(\epsilon^\prime)^2$ as: \[x^\prime\mapsto \alpha_1\epsilon^\prime,\quad y^\prime\mapsto \beta_1\epsilon^\prime,\quad z^\prime\mapsto \epsilon^\prime.\]
Generators of the fibre product are given by $x=x^\prime-\alpha_1z^\prime,\ y=y^\prime-\beta_1z^\prime$ as above. %- these suffice, since $z^\prime\m^\prime\subseteq\langle x^\prime,y^\prime\rangle^2$ and $\langle x^\prime,y^\prime\rangle^2=\langle x,y\rangle^2$ (using completeness).
 \end{remark}
 
 \begin{remark}
 If $A_3(1,1)=0$ (so also $A_4(2,1)=\Delta=0$), the resulting singularity is hyperelliptic. Indeed, consider the involution $\iota$ of $\CC[\![t_1]\!]\oplus\ldots\oplus\CC[\![t_b]\!]$ defined by:
  \[\left\{\begin{aligned}
            t_4 &\leftrightarrow A_4(1,1)t_1+A_4(1,2)t_1^2  &\\
            t_3  &\leftrightarrow A_3(2,1)t_2-\frac{A_3(1,2)}{B_2(1,2)}t_2^2&\\
            t_i&\leftrightarrow -\frac{A_i(1,2)}{A_{i+1}(1,2)}t_{i+1} & \qquad i=5,7,\ldots,b-1
           \end{aligned}\right.
 \]
 if $b$ is even; otherwise extend to $b$ by setting $t_{b+1}=t_b$ and $A_{b+1}=A_b$ formally, so to send $t_b$ to $-t_b$. This defines an automorphism of the normalisation since every term in the expression above is a parameter for the corresponding branch. Observe then that $\iota$ preserves the subalgebra $R$ and makes the differential $\eta$ anti-invariant. %We conclude that this family is associated to the hyperelliptic component of the stratum $\Hcal(1,1,1,1)$.
 
 We may also check that these singularities are hyperelliptic by matching the local generators with the expressions given in \cite[\S3]{BB23}. For simplicity, set all crimping parameters to $1$ or $0$ as appropriate. Consider the change of coordinates:
 
  \[\left\{\begin{aligned}
            s_1 &=A_4                     &=&t_1+t_4,                               &\\
            s_2 &=A_3 &=&t_2+t_3, \\
            u_2 &= 2B_2-A_3^2-A_4^2 &=&t_1^2+t_2^2-t_3^2-t_4^2, \\
            s_i &=A_{2i}-A_{2i+1}         &=&t_{2i}-t_{2i+1},                       &\\
            u_i &=A_{2i}+A_{2i+1}-A_4^2   &=&t_1^2-t_4^2+t_{2i}+t_{2i+1},           & i=3,4,\ldots,\lfloor \frac{b}{2}\rfloor\\
           \end{aligned}\right.
 \]
  if $b$ is even, and adjoin 
  \[\left\{\begin{aligned}
            s_{b'} &=A_b^2-A_4^4                     &=&\ -t_4^4+t_b^2,     \\
            u_{b'} &=A_b   &=&t_1^2+t_b,           \\
           \end{aligned}\right.
 \]
 in case $b$ is odd (where $b'=\frac{b+1}{2}$).
 The following vanish: $u_2^2-s_1^4-s_2^4,u_i^2-s_1^4-s_i^2,s_1(u_i-u_j),s_iu_j$, $u_iu_j=s_1^4$, and $u_{b'}^2-s_1^4-s_{b'}$ in case $b$ is odd. The corresponding tropical differential has slope $2$ along the edges $e_1,e_4$ (conjugate), $e_2,e_3$ (conjugate), and slope $1$ along all the other edges.
 
 Note that for $b=4$ the embedding dimension is $3$ in this case.
\end{remark}

\subsection*{[2 0 1]} In this case $\tm^4\subseteq R$ and $b\geq 2$.
Let us assume that $t_1^3$ does not belong to $R$, then neither does $t_1$; another linear monomial is missing, which we assume to be $t_2$. We can then find generators:
\begin{align*}
 C_i=&t_i^3+C_i(1,3)t_1^3,&i=2,\ldots,b;\\
 B_i=&t_i^2+B_i(1,3)t_1^3,&i=1,\ldots,b;\\
 A_i=&t_i+A_i(1,1)t_1+A_i(1,3)t_1^3+A_i(2,1)t_2,&i=3,\ldots,b.
\end{align*}

Observe that $A_iB_1=A_i(1,1)t_1^3$, hence $A_i(1,1)=0$ for all $i$. Then $A_iB_2=A_i(2,1)t_2^3$ must be a multiple of $C_2$, so $C_2(1,3)=0$ or $A_i(2,1)=0$ for all $i$.
\begin{itemize}[leftmargin=0cm]
 \item If $C_2(1,3)=0$, then \[A_i^3=C_i+A_i(2,1)^3C_2\Leftrightarrow C_i(1,3)=0.\]
 So $C_i=t_i^3$ for all $i=2,\ldots,b$.
 
 Suppose that $A_i(2,1)=0$ for all $i=3,\ldots,b$. Then, $B_i-A_i^2=B_i(1,3)t_1^3\pmod{\tm^4}$ implies $B_i(1,3)=0$ for all $i=3,\ldots,b$. Hence $R/\cond$ would be generated by $1,A_3,\ldots,A_b,B_1,B_2$ and of dimension $\leq b+1$, which is too small.
 
 We may therefore assume $b\geq3$ and $A_3(2,1)\neq0$. If there was another index $j$ (say $j=4$) such that $A_4(2,1)\neq 0$, then $A_3A_4$ would be a non-zero multiple of $t_2^2$. Subtracting a multiple of $t_2^2$ from $A_i^2$ we would find $t_i^2\in R$ for all $i=3,\ldots,b$ as well. Hence $R/\cond=\langle 1,A_3,\ldots,A_b,B_1\rangle$, which is again too small.
 Summing up, the generators of $\m/\m^2$ up to $\tm^4$ are:
\begin{align*}
  A_3= && A_3(1,3)t_1^3 && A_3(2,1)t_2 && t_3 && 0 && \ldots && 0 \\
  A_4= && A_4(1,3)t_1^3 && 0 && 0 && t_4 && \ldots && 0 \\
  \ldots &&&&&&&&&& \\
  A_b= && A_b(1,3)t_1^3 && 0 && 0 && 0 && \ldots && t_b \\
  B_1= && t_1^2+B_1(1,3)t_1^3 && 0 && 0 && 0 && \ldots && 0 \\
  B_2= && B_2(1,3) t_1^3 && t_2^2 && 0 && 0 && \ldots && 0 \\
 \end{align*}
 with $A_3(2,1), B_2(1,3)$ and $A_i(1,3),\ i=4,\ldots,b,$ invertible (by indecomposability) and $A_3(1,3),B_1(1,3)$ arbitrary. The crimping space is isomorphic to $\Gm^{b-1}\times \Ga^2$. We verify that $R/\cond=\langle1,A_3,\ldots,A_b,B_1,B_2,A_3^2\rangle.$

 \begin{comment}
 \begin{rem}
  For $b=3$, this singularity consists of a cusp (first branch) and two smooth branches tangent to each other (a tacnode) glued together non-transversely (so not in $\Aaff^{2+2}$ but in $\Aaff^3$). Setting all crimping parameteres to $0$ or $1$, and $x=A_3,y=B_1,z=B_2$, local equations are $(xy,y^3+x^2z-z^2)$.
 \end{rem}
 \end{comment}

 The dualising line bundle is generated by the meromorphic differential:
\[(1-B_1(1,3)t_1)\frac{\de t_1}{t_1^4}-B_2(1,3)\frac{\de t_2}{t_2^3}+B_2(1,3)A_3(2,1)^2\frac{\de t_3}{t_3^3}-\sum_{i=3}^bA_i(1,3)\frac{\de t_i}{t_i^2}.\]

It follows that in the semistable tail:
\[\omega_Z=\OO_Z(2p_1+p_2+p_3),\]
so this singularity corresponds to the stratum $\Hcal(2,1,1)$.

\begin{remark}
 In case $B_1(1,3)=0$, these singularities are hyperelliptic. Indeed, consider the involution $\iota$ of $\CC[\![t_1]\!]\oplus\ldots\oplus\CC[\![t_b]\!]$ defined by:
  \[\left\{\begin{aligned}
            t_1 &\leftrightarrow -t_1  &\\
            t_3  &\leftrightarrow A_3(2,1)t_2-\frac{A_3(1,2)}{B_2(1,3)}t_2^2&\\
            t_i&\leftrightarrow -\frac{A_i(1,2)}{A_{i+1}(1,2)}t_{i+1} & \qquad i=5,7,\ldots,b-1
           \end{aligned}\right.
 \]
 if $b$ is even; otherwise extend to $b$ by setting $t_{b+1}=t_b$ and $A_{b+1}=A_b$ formally, so to send $t_b$ to $-t_b$. This defines an automorphism of the normalisation since every term in the expression above is a parameter for the corresponding branch. Observe then that $\iota$ preserves the subalgebra $R$ and makes the differential $\eta$ anti-invariant.
 
 The corresponding tropical differential has slope $3$ along the edge $e_1$ (Weierstrass), $2$ along $e_2,e_3$ (conjugate), and slope $1$ along all the other edges.
\end{remark}

\begin{remark}
 We sketch an alternative proof that this is the only family of singularities corresponding to the stratum $\Hcal(2,1,1)$. We want $\omega_{C,p}$ to be generated by a differential of the form:
 \[\eta=\frac{\de t_1}{t_1^4}(1+a_1t_1+a_2t_1^2)+\frac{\de t_2}{t_2^3}(u+b_1t_2)+\frac{\de t_3}{t_3^3}(v+c_1t_3)+\sum_{i=4}^bu_i\frac{\de t_i}{t_i^2},\]
 with $u,v,u_i$ invertible.
 We know then that \[\cond=\langle t_1^4,t_2^3,t_3^3,t_4^2,\ldots,t_b^2\rangle.\]
 
 Let $S$ denote the local ring of the seminormalisation, $S=\CC[\![t_1,\ldots,t_b]\!]/(t_it_j)$. The following diagram is Cartesian and displays the singularity as the result of pinching the seminormalisation along a finite subscheme:
 \bcd
 R\ar[r,hook]\ar[d]\ar[dr,phantom,"\ulcorner"] & S\ar[d]\\
 R/\cond\ar[r,hook] & S/\cond.
 \ecd
 Therefore, we have to look for subalgebras of \[S/\cond=\CC[\![t_1,\ldots,t_b]\!]/(t_1^4,t_2^3,t_3^3,t_4^2,\ldots,t_b^2,t_it_j)\] of codimension $g$. Moreover, they have to be isotropic with respect to the residue pairing $\langle f,g\rangle=\on{Res}(fg\ \eta)$.
 
 Suppose that we had a generator $L_1$ linear in $t_1$. Then $L_1^3=t_1^3$ would have residue $1$ with $\eta$. For the same reason, $t_2$ and $t_3$ may only appear in combination $vt_2-ut_3$. Finally, $t_1^3$ may not appear: this exhausts the $3$-dimensional gap. Notice that the correct power of $t_1^3$ must be appended to every generator in order to ensure that the residue vanishes. We recover the same generators as above.
\end{remark}

\medskip

\item Assume instead $C_2(1,3)\neq 0$ and $A_i(2,1)=0$ for all $i$. In this case, the generators of $\m/\m^2$ up to $\tm^4$ are:

\begin{align*}
  A_3= && A_3(1,3)t_1^3 && 0 && t_3  && \ldots && 0 \\
  \ldots &&&&&&&&&& \\
  A_b= && A_b(1,3)t_1^3 && 0 &&  0 && \ldots && t_b \\
  B_1= && t_1^2+B_1(1,3)t_1^3 && 0 && 0 &&  \ldots && 0 \\
  B_2= && B_2(1,3) t_1^3 && t_2^2 && 0 && \ldots && 0 \\
  C_2= && C_2(1,3) t_1^3 && t_2^3 && 0 && \ldots && 0 \\
 \end{align*}
 with $C_2(1,3)$ and $A_i(1,3)$ invertible for $i=3,\ldots,b$, and $B_1(1,3),B_2(1,3)$ arbitrary. The crimping space is isomorphic to $\Gm^{b-1}\times \Ga^2$. We verify that $R/\cond=\langle1,A_3,\ldots,A_b,B_1,B_2,C_2\rangle.$
 
 \begin{rem}
  For $b=2$, the singularity consists of two cusps glued together non-transversely in $\Aaff^3$.
 \end{rem}

 The dualising line bundle is generated by the meromorphic differential:
\[(1-B_1(1,3)t_1)\frac{\de t_1}{t_1^4}-(C_2(1,3)+B_2(1,3)t_2)\frac{\de t_2}{t_2^4}-\sum_{i=3}^bA_i(1,3)\frac{\de t_i}{t_i^2}.\]

It follows that in the semistable tail:
\[\omega_Z=\OO_Z(2p_1+2p_2).\]
 
 \begin{lem}\label{lem:theta_odd}
  $p_1+p_2$ is an \emph{odd} theta characteristic on $Z$.
 \end{lem}

 We consider a rational compactification of $(C,p)$, i.e. a projective curve $C$ of arithmetic genus three having $p$ as its only singularity. We study theta characteristics on $C$ as in \cite{HarrisTheta}. We refer to \S\ref{sec:Harris} for the notation. Since the adjoint divisor is even, the linear form $l$ vanishes identically. Let $I'$ be the ideal of $\OO_C$ that pulls back to $\OO(-p_1-p_2)$ on the normalisation. In our case $I'=\m$. The number of even and odd theta characteristics is determined by $\on{len}(\OO_C/I')=1$. In particular, there is exactly one odd theta characteristic $\theta$ (the sum of the two points at infinity on the cuspidal branches), and no even ones. Consider a smoothing of $(C,\theta)$. By semistable reduction, we can substitute the central fibre with a nodal curve $C^{ss}$; we may assume that there exists a contraction $C^{ss}\to C$. By collapsing the rational tails of $C^{ss}$, we see $p_1+p_2$ as the limit of $\theta$ on the genus three semistable tail $Z$.

 \smallskip
 
Summing up, this family corresponds to $\Hcal(2,2)^{\text{odd}}$.
\end{itemize}

\subsection*{[1 1 1]} In this case too $\tm^4\subseteq R$ and $b\geq 2$. Let us say that $t_1^3$ is missing; then so is $t_1$. If the missing quadratic term is $t_1^2$, then we can write our linear generators as $A_i=t_i+\sum_{j=1}^3A_i(1,j)t_1^j$ for $i=2,\ldots,b$. Suppose instead that we do not have $t_2^2$. We have though a linear generator of the form $L_2=t_2+L_2(1,1)t_1+h.o.t.$ with $L_2(1,1)\neq0$ (otherwise $L_2^2=t_2^2+h.o.t.$); but then we can add a multiple of $L_2^2$ to any generator of $R$ in order to change their $t_2^2$ term into a $t_1^2$ term. Hence we can indeed write the linear generators $A_i$ as above.

If $b\geq3$, by pairing $A_i$ with $A_j$ for $i\neq j$ we see that there is at most one $i$ such that $A_i(1,1)\neq 0$. Indeed, there must be exactly one, otherwise $B_i=t_i^2$ and $C_i=t_i^3$, so $R/\cond=\langle 1,A_2,\ldots,A_b\rangle$ would be too small.

Ultimately, the generators of $\m/\m^2$ are:
 \begin{align*}
  A_2= && A_2(1,1)t_1+A_2(1,2)t_1^2+A_2(1,3)t_1^3 && t_2 && 0 && \ldots && 0 \\
  A_3= && A_3(1,3)t_1^3 && 0 && t_3 && \ldots && 0 \\ 
  \ldots &&&&&&&&&& \\
  A_b= && A_b(1,3)t_1^3 && 0  && 0 && \ldots && t_b \\
 \end{align*}
with $A_2(1,1)$ and $A_i(1,3)$ invertible for $i=3,\ldots,b$, and $A_2(1,2),A_2(1,3)$ arbitrary. The crimping space is isomorphic to $\Gm^{b-1}\times \Ga^2$. We verify that $R/\cond=\langle1,A_2,\ldots,A_b,A_2^2,A_2^3\rangle.$
 
 \begin{rem}\label{rem:hyperelliptic}
  For $b=2$, we need one more generator, for instance $D=t_1^4-t_2^4\in\tm^4$. This is the $A_7$-singularity $\CC[\![A,D]\!]/(D^2-A^8)$. It is hyperelliptic.
  
  More generally, these singularities are hyperelliptic in the sense of \cite{BB23}. In order to see it, set all crimping parameters to either $1$ or $0$ for simplicity. The corresponding tropical differential has slope $3$ on the branches corresponding to $t_1$ and $t_2$, and $1$ on all other branches: if we have an even number of branches $b=2b'$, then we can assume that all branches come in conjugate pairs; otherwise, we assume that the branch corresponding to $t_{2b+1}$ is Weierstrass. In the former case we can choose new coordinates:
  \[\left\{\begin{aligned}
            s_1 &=A_2                     &=&t_1+t_2,                               &\\
            s_i &=A_{2i-1}-A_{2i}         &=&t_{2i-1}-t_{2i},                       &\\
            u_i &=A_{2i-1}+A_{2i}-A_2^3   &=&t_1^3-t_2^3+t_{2i-1}+t_{2i},           & i=2,\ldots,b',\\
           \end{aligned}\right.
 \]
  satisfying $u_i^2-s_1^6-s_i^2,s_1(u_i-u_j),s_iu_j$, and $u_iu_j=s_1^6$, see \cite[\S 3]{BB23}. In the latter case, we adjoin 
  \[\left\{\begin{aligned}
            s_{b'+1} &=A_b^2                     &=&\ t_1^6+t_b^2,     \\
            u_{b'+1} &=A_b-A_2^3   &=&-t_2^3+t_b,           \\
           \end{aligned}\right.
 \]
 satisfying $u_{b'+1}^2+s_1^6-s_{b'+1}$ and all the other equations as above.
  \end{rem}

The dualising line bundle is generated by the meromorphic differential $\eta$:

\[\left(1-2\frac{{A_2(1,2)}}{{A_2(1,1)}}t_1-\left(\frac{{A_2(1,3)}}{{A_2(1,1)}}-2\frac{{A_2(1,2)}^2}{{A_2(1,1)}^2}\right)t_1^2\right)\frac{\de t_1}{t_1^4}-{A_2(1,1)}^3\frac{\de t_2}{t_2^4}-\sum_{i=3}^bA_i(1,3)\frac{\de t_i}{t_i^2}.\]

It follows that in the semistable tail the attaching points of the two osculating branches on the core $Z$ are such that:
\[\omega_Z=\OO_Z(2p_1+2p_2).\]
\begin{lem}\label{lem:even_theta}
 $p_1+p_2$ is an even theta characteristic on $Z$.
\end{lem}
 Following Harris, we see once again that $l\equiv 0$, but in this case $R/I^\prime=\CC\langle 1,x\rangle$, so there is exactly one even theta characteristic $\theta$, and no odd ones on $C$. We conclude by semistable reduction of $(C,\theta)$ as above.

 We offer an alternative argument. Let $C$ be the genus three compactification as above. Reasoning as in \S\ref{sec:sstails}, we observe that $\{\eta,A_2\eta,A_2^2\eta\}$ extends to a basis of $H^0(C,\omega_C)$. Applying semistable reduction to $A_2\eta$ (resp. $A_2^2\eta$), we find a differential on the semistable tail $Z$ with multiplicity $(1,1)$ (resp. $(0,0)$) at $(p_1,p_2)$. So, $p_1$ and $p_2$ do not impose independent conditions on $H^0(\omega_Z)$. This is the case when $p_1+p_2$ is an even theta characteristic (two conjugate points on a hyperelliptic curve), but not when it is odd (bitangent of a quartic, or, specialising, sum of two distinct Weierstrass points on a hyperelliptic curve). See Table \ref{rkfun} for the corresponding rank functions.
 
 \begin{table}[h]
 \begin{tabular}{c|ccc}
 $H^0(\omega_Z(-))$ & $0$ & $-p_1$ & $-2p_1$ \\
 \hline
  $0$     & $3$ & $2$ & $1$ \\
  $-p_2$  & $2$ & $2$ & $1$ \\
  $-2p_2$ & $1$ & $1$ & $1$ \\
 
 \end{tabular}\qquad
 \begin{tabular}{c|ccc}
 $H^0(\omega_Z(-))$ & $0$ & $-p_1$ & $-2p_1$ \\
 \hline
  $0$     & $3$ & $2$ & $1$ \\
  $-p_2$  & $2$ & $1$ & $1$ \\
  $-2p_2$ & $1$ & $1$ & $1$ \\
 \end{tabular}
 \medskip
 \caption{Left: rank function of an even theta. Right: odd theta.}\label{rkfun}
 \end{table}
Summing up, these singularities correspond to $\Hcal(2,2)^{\text{even}}$.

\subsection*{[1 1 0 1]} In this case $\tm^5\subseteq R$. The monomial unibranch singularity $\CC[\![t^3,t^5,t^7]\!]$ is not Gorenstein: for instance, the associated monoid is not symmetric. So assume $b\geq2$. The arguments are similar to the previous case: if we miss the monomial $t_1^4$, then we also miss $t_1$ and $t_1^2$. We have linear generators $A_i=t_i+\sum_{j=1,2,4}A_i(1,j)t_1^j,\ i=2,\ldots,b$, and a cubic generator $C_1=t_1^3+C_1(1,4)t_1^4$. In order for $A_iC_1$ not to start with $t_1^4$, we need $A_i(1,1)=0$ for all $i$. For the same reason, there is at most one index $i$ such that $A_i(1,2)\neq 0$. In fact there is exactly one, otherwise $R/\cond=\langle 1,A_2,\ldots,A_b,C_1\rangle$, which is too small. So we have generators:
\begin{align*}
  A_2= && A_2(1,2)t_1^2+A_2(1,4)t_1^4 && t_2 && 0 && \ldots && 0 \\
  A_3= && A_3(1,4)t_1^4 && 0 && t_3 && \ldots && 0 \\
   \ldots &&&&&&&&&& \\
  A_b= && A_b(1,4)t_1^4 && 0 && 0 && \ldots && t_b \\
  C_1= && t_1^3+C_1(1,4)t_1^4 && 0 && 0 && \ldots && 0 \\
\end{align*}
 with $A_2(1,2)$ and $A_i(1,4),\ i=3,\ldots,b,$ invertible, and $A_2(1,4),C_1(1,4)$ arbitrary. The crimping space is isomorphic to $\Gm^{b-1}\times \Ga^2$. We verify that $R/\cond=\langle1,A_2,\ldots,A_b,A_2^2,C_1\rangle.$
 
 \begin{rem}
  For $b=2$, setting $x=A_2,y=C_1$, we find local equations $y(y^2-x^3)$. This is the $E_7$-singularity.
 \end{rem}

The dualising line bundle is generated by the meromorphic differential:
\[\left(1-C_1(1,4)t_1-\frac{A_2(1,4)}{A_2(1,2)}t_1^2\right)\frac{\de t_1}{t_1^5}
-A_2(1,2)^2\frac{\de t_2}{t_2^3}
-\sum_{i=3}^bA_i(1,4)\frac{\de t_i}{t_i^2}.
\]

We can thus write:
\[\omega_Z=\OO_Z(3p_1+p_2),\]
so this singularity corresponds to the stratum $\Hcal(3,1)$.

\subsection*{[1 1 0 0 1]} In this case $\tm^6\subseteq R$. We encounter the monomial unibranch singularity $\CC[\![t^3,t^4]\!]$, classically known as $E_6$. In general, we have linear generators $A_i=t_i+\sum_{j=1,2,5}A_i(1,j)t_1^j,\ i=2,\ldots,b$, but also a cubic $C_1=t_1^3+C_1(1,5)t_1^5$ and a quartic $D_1=t_1^4+D_1(1,5)t_1^5$. Multiplying $A_i$ with $D_1$ and $C_1$, we see that both $A_i(1,1)$ and $A_i(1,2)$ must vanish. So:
\begin{align*}
 A_2= && A_2(1,5)t_1^5 && t_2 && 0 && \ldots && 0 \\
   \ldots &&&&&&&&&& \\
 A_b= && A_b(1,5)t_1^5 && 0 && 0 && \ldots && t_b \\
 C_1= && t_1^3+C_1(1,5)t_1^5 && 0 && 0 && \ldots && 0 \\
 D_1= && t_1^4+D_1(1,5)t_1^5 && 0 && 0 && \ldots && 0 \\
 \end{align*}
 with $A_i(1,5)\in\Gm,\ i=2,\ldots,b,$ and $C_1(1,5),D_1(1,5)\in\Ga$. We check:
$R/\cond=\langle1,A_2,\ldots,A_b,C_1,C_2\rangle.$ 

The dualising line bundle is generated by the meromorphic differential:
\[\left(1-D_1(1,5)t_1-C_1(1,5)t_1^2\right)\frac{\de t_1}{t_1^6}
-\sum_{i=2}^b A_i(1,5)\frac{\de t_i}{t_i^2}.
\]

It follows that we can write:
\[\omega_Z=\OO_Z(4p_1).\]
\begin{lem}
 $2p$ is an odd theta characteristic.
\end{lem}
 Similar to Lemma \ref{lem:theta_odd}. So this singularity corresponds to $\Hcal(4)^{\text{odd}}$.

\subsection*{[1 0 1 0 1]} In this case $\tm^6\subseteq R$. We have a monomial unibranch singularity $\CC[\![t^2,t^7]\!]$, classically known as $A_6$, which is hyperelliptic. For $b\geq2$, there are linear generators $A_i=t_i+\sum_{j=1,3,5}A_i(1,j)t_1^j,\ i=2,\ldots,b$, and a quadratic one $B_1=t_1^2+B_1(1,3)t_1^3+B_1(1,5)t_1^5$. Multiplying $A_i$ with $B_1^2$ and with $B_1$, we see that both $A_i(1,1)$ and $A_i(1,3)$ must vanish. So:
\begin{align*}
 A_2= && A_2(1,5)t_1^5 && t_2 && 0 && \ldots && 0 \\
   \ldots &&&&&&&&&& \\
 A_b= && A_b(1,5)t_1^5 && 0 && 0 && \ldots && t_b \\
 B_1= && t_1^2+B_1(1,3)t_1^3+B_1(1,5)t_1^5 && 0 && 0 && \ldots && 0 \\
 \end{align*}
 with $A_i(1,5)\in\Gm,\ i=2,\ldots,b,$ and $B_1(1,3),B_1(1,5)\in\Ga$. We check:
$R/\cond=\langle1,A_2,\ldots,A_b,B_1,B_1^2\rangle.$

The dualising line bundle is generated by the meromorphic differential:
\[\left(1-2B_1(1,3)t_1-B_1(1,5)t_1^3\right)\frac{\de t_1}{t_1^6}
-\sum_{i=2}^b A_i(1,5)\frac{\de t_i}{t_i^2}.
\]

It follows that we can write:
\[\omega_Z=\OO_Z(4p).\]
\begin{lem}
 $2p$ is an even theta characteristic.
\end{lem}
Similar to Lemma \ref{lem:even_theta}. So this singularity corresponds to $\Hcal(4)^{\text{even}}$.

\begin{rem}
 These singularities are hyperelliptic in the sense of \cite{BB23}. For simplicity, set all crimping parameters equal to $1$ or $0$, and assume that $b=2b'+1$. The associated tropical differential has slope $5$ along the edge corresponding to $p_1$, which is Weierstrass. With coordinates:
 \[\left\{\begin{aligned}
            s_1 &=B_1                                &=&\ t_1^2,                               &\\
            s_{i+1} &=\frac{1}{2}(A_{2i+1}-A_{2i})       &=&\ \frac{1}{2}(t_{2i+1}-t_{2i}),                       &\\
            u_{i+1} &=\frac{1}{2}(A_{2i+1}+A_{2i})       &=&\ t_1^5+\frac{1}{2}(t_{2i+1}+t_{2i}),           & i=1,\ldots,b',\\
           \end{aligned}\right.
 \]
  we find the local equations $u_i^2-s_1^5-s_i^2,s_1(u_i-u_j),s_iu_j$, and $u_iu_j=s_1^6$, see \cite[\S 3]{BB23}. The case $b$ even can be addressed by adding another Weierstrass branch as in Remark \ref{rem:hyperelliptic}.
\end{rem}

\subsection*{Funding}
 During the preparation of this work, I was partially supported by the Deutsche Forschungsgemeinschaft (DFG, German Research Foundation) under Germany’s Excellence Strategy EXC-2181/1 - 390900948 (the Heidelberg STRUCTURES Cluster of Excellence), and by the European Union - NextGenerationEU under the National Recovery and Resilience Plan (PNRR) - Mission 4 Education and research - Component 2 From research to business - Investment 1.1 Notice Prin 2022 - DD N. 104 del 2/2/2022, from title "Symplectic varieties: their interplay with Fano manifolds and derived categories", proposal code 2022PEKYBJ – CUP J53D23003840006.

\subsection*{Conflict of interests} The author has no conflicts of interest to declare that are relevant to the content of this article.
 
%\bibliographystyle{alpha}
%\bibliography{biblio} 

%\begin{comment}

%\end{comment}

 \bigskip
 
\end{document}